%% file: Humeri.tex
\documentclass[ DIV=15, font size=10.5pt,
                 onecolumn,
               ]{scrartcl}  

\input{utils/packages}
\input{utils/pgfplotsettings}

\input{utils/commonCommands}

\input{utils/titleAndAuthors}
\input{utils/journalFooter}

\crefname{figure}{Fig.}{Fig.}
\crefname{equation}{Eq.}{Eq.}
\crefname{table}{Tab.}{Tab.}

\drawboundingboxtrue
\drawboundingboxfalse 

\newcommand*{\figref}[2][]{%
	\hyperref[{fig:#2}]{%
		Fig.~\ref*{fig:#2}%
		\ifx\\#1\\%
		\else
		\,#1%
		\fi
	}%
}

\definecolor{changes}{RGB}{0,0,0}
\definecolor{reviewer}{RGB}{0,0,255}
\definecolor{reviewer2}{rgb}{0,0.6,0}
\definecolor{todos}{RGB}{255,0,0}

\newfloatcommand{capbtabbox}{table}[][\FBwidth]

\begin{document}  
\normalem
\maketitle  
  
%% Abstract ---------------------------------------
\vspace{-1.5cm} 
\hrule 
\section*{Abstract}
\input{utils/abstract}
 \vspace{.2cm} 
%% Keywords ---------------------------------------
\vspace{0.25cm}\\
\noindent \textit{Keywords:} \input{utils/keywords} 
\vspace{0.35cm}
\hrule 
\vspace{0.15cm}
\captionsetup[figure]{labelfont={bf},name={Fig.},labelsep=colon}
\captionsetup[table]{labelfont={bf},name={Tab.},labelsep=colon}
\tableofcontents
\vspace{0.5cm}
\hrule 

%% Actual Content ---------------------------------
\input{sections/introduction}
\input{sections/methods}
\input{sections/results}

\input{sections/conclusions}

%% Acknowledgements -------------------------------  
\newpage
\section*{Acknowledgements} 
\input{utils/acknowledgements}

\section*{Conflict of interest} 
No potential conflict of interest was reported by the authors.

\newpage
%% References -------------------------------------
\bibliographystyle{apalike}
\bibliography{library.bib}
 
\newpage
\input{utils/appendix}

\end{document}

%% file: utils/packages.tex
\usepackage{color} 
\usepackage{psfrag} 
\usepackage{caption}
\usepackage{hyperref}
\usepackage{amsmath}
\usepackage{amssymb}
\usepackage{graphics}
\usepackage{mathtools}
 \usepackage{pstool} 
\usepackage{todonotes}
\usepackage{gensymb}
\usepackage{bm}
\usepackage{floatrow}
\usepackage{multirow}
\usepackage{dsfont}

\usepackage{textcomp}

\usepackage[square, numbers,sort&compress]{natbib}
 
\usepackage{pgfplots}
\usepackage{pgfplotstable}
\usepackage{filecontents}

\usepackage{tikz}
\usepackage{tikzscale}

\usetikzlibrary{spy}
\usetikzlibrary{snakes,arrows,shapes}    
\usetikzlibrary{spy}
\usetikzlibrary{backgrounds}  
\usetikzlibrary{decorations}
\usetikzlibrary{positioning}
\usetikzlibrary{patterns}
\usetikzlibrary{calc} 
\usepackage{booktabs}   
\usepackage{makecell} 
\usepackage{colortbl}   

\usepackage{xspace}
\usepackage{color}     
\usepackage{setspace}   
 
\usepackage{soul}
\usepackage{ulem}
\usepackage{currfile}
\usepackage{import}

\usepackage{authoraftertitle} % used for getting the authors 
\usepackage{authblk} % used for affiliations 
\usepackage{cleveref}

%% file: utils/pgfplotsettings.tex
\pgfplotsset{compat=1.8}

\newcommand{\findmax}[3]{
    \pgfplotstablesort[sort key={#2},sort cmp={float >}]{\sorted}{#1}%
    \pgfplotstablegetelem{0}{#2}\of{\sorted}%
    \let #3=\pgfplotsretval%
}

\pgfplotscreateplotcyclelist{myCycleListBlackAndWhite}
{%
  {black, mark=o},
  {black, mark=square},
  {black, mark=triangle},
  {black, mark=diamond}, 
  {black, mark=pentagon}, 
  {black, mark=*},
  {black, mark=square*},
  {black, mark=triangle*},
  {black, mark=diamond*}, 
  {black, mark=pentagon*}, 
}

\definecolor{darkgreen}{rgb}{0,0.4,0} 
\definecolor{darkbrown}{rgb}{0.5, 0.396, 0.09}

\definecolor{c1}{rgb}{0.0, 0.4196078431372549, 0.6431372549019608}
\definecolor{c2}{rgb}{1.0, 0.5019607843137255, 0.054901960784313725}
\definecolor{c3}{rgb}{0.6705882352941176, 0.6705882352941176,
0.6705882352941176} \definecolor{c}{rgb}{0.34901960784313724, 0.34901960784313724, 0.34901960784313724}
\definecolor{c4}{rgb}{0.37254901960784315, 0.6196078431372549,
0.8196078431372549} \definecolor{c}{rgb}{0.7843137254901961, 0.3215686274509804, 0.0}
\definecolor{c5}{rgb}{0.5372549019607843, 0.5372549019607843,
0.5372549019607843} \definecolor{c}{rgb}{0.6352941176470588, 0.7843137254901961, 0.9254901960784314}
\definecolor{c6}{rgb}{1.0, 0.7372549019607844, 0.4745098039215686}
\definecolor{c7}{rgb}{0.8117647058823529, 0.8117647058823529,
0.8117647058823529}

\pgfplotscreateplotcyclelist{myCycleListColor}
{%
  {blue, mark=o},
  {red, mark=square},
  {darkbrown, mark=triangle},
  {darkgreen, mark=diamond},
  {black, mark=pentagon},
  {blue, mark=*},
  {red, mark=square*},
  {darkbrown, mark=triangle*},
  {darkgreen, mark=diamond*},
  {black, mark=pentagon*},
}

\pgfplotsset{every axis/.append style= 
              {
                font=\small,
                mark size=2,
                line width = 0.1,
                legend style={font=\small, mark size=3, draw=none, fill=none},
                legend cell align=left,
                cycle list name=myCycleListColor,
              }
            }
            
\pgfplotstableset
{
  every odd row/.style={before row={\rowcolor[gray]{0.9}}},
  every head row/.style=
  {
    before row=
    {%
      \toprule
    },
    after row=\midrule
  },
  every last row/.style=
  {
    after row=\bottomrule
  }
}

\newif\ifdrawboundingbox

\usetikzlibrary{external} % Is activated with the command
\tikzset{external/system call={pdflatex \tikzexternalcheckshellescape
-halt-on-error -interaction=batchmode -jobname "\image" "\texsource"}} 
\pgfkeys
{ 
  /pgf/number format/.cd,
  fixed,
  set thousands separator={\,},
  1000 sep in fractionals,
}

%% file: utils/commonCommands.tex
\newcolumntype{C}[1]{>{\centering\arraybackslash}m{#1}}
\newcolumntype{R}[1]{>{\raggedright\arraybackslash}m{#1}}
\newcolumntype{L}[1]{>{\raggedleft\arraybackslash}m{#1}}

\newcommand{\delete}[1]{\xspace}

%% file: utils/titleAndAuthors.tex
\title{Predicting Fracture in the Proximal Humerus using Phase Field Models}
\author[1]{L. Hug\thanks{lisa.hug@tum.de, Corresponding Author}}
\author[2]{G. Dahan}%\thanks{second.author@xyz.de}}
\author[1]{S. Kollmannsberger}%\thanks{second.author@xyz.de}}
\author[1]{E. Rank}%\thanks{second.author@xyz.de}}
\author[2]{Z. Yosibash}%\thanks{second.author@xyz.de}}

 \affil[1]{Chair for Computational Modeling and Simulation,
 Technical University of Munich, Arcisstr. 21, 80333 Munich, Germany}
 \affil[2]{School of Mechanical Engineering, The Iby and Aladar Fleischman Faculty of Engineering, Tel-Aviv University, 69978 Ramat-Aviv, Israel}

\newcommand{\publicationDate}{\today}

%% file: utils/journalFooter.tex
\usepackage{fancyhdr}
\date{}
\fancypagestyle{plain}{%
  \fancyhf{}%
  \fancyfoot[L]{
  \vspace{-1.25cm} 
  \footnotesize{Preprint submitted to \textit{Journal}  \hfill
  \publicationDate
  \\
  }} }

%% file: utils/abstract.tex
Proximal humerus impacted fractures are of clinical concern in the elderly population. Prediction of such fractures by CT-based finite element methods encounters several major obstacles such as heterogeneous mechanical properties and fracture due to compressive strains. We herein propose to investigate a variation of the phase field method (PFM) embedded into the finite cell method (FCM) to simulate impacted humeral fractures in fresh frozen human humeri. The force-strain response, failure loads \textcolor{changes}{and the fracture path} are compared to experimental observations for validation purposes. \textcolor{changes}{The PFM (by means of the regularization parameter $\ell_0$)} is first calibrated by one experiment and thereafter used for the prediction of the mechanical response
of two other human fresh frozen humeri. All humeri are fractured at the surgical neck and strains are monitored by Digital Image Correlation~(DIC). \\ Experimental strains in the elastic regime are reproduced with good agreement ($R^2 = 0.726$), similarly to the validated finite element method \cite{dahan2020neck}. The failure pattern and fracture evolution at the surgical neck predicted by the PFM mimic extremely well the experimental observations for all three humeri. The maximum relative error in the computed failure loads \textcolor{changes}{is $3.8\%$}. To the best of our knowledge this is the first method that can predict well the experimental compressive failure pattern as well as the force-strain relationship \textcolor{changes}{in proximal humerus fractures}. \\[1em]

%% file: utils/keywords.tex
humerus fracture, brittle fracture, phase-field modeling, Finite Cell Method, iso-geometric analysis

%% file: sections/introduction.tex
\section{Introduction} \label{sec:intro}
Osteoporotic fractures are a frequent injury among the elderly and often require medical attention. \textcolor{changes}{With an incidence of $5-6\%$ of all adult fractures \cite{court2012epidemiology}, proximal humerus fractures are one of the most common type of fragility fractures in patients over 65, next to the proximal femur, vertebral body and distal radius  \cite{kim2012epidemiology, tsuda2017epidemiology}.} Proximal humerus fractures are often caused by low-impact falls on an outstretched arm, which induce a fracture at the surgical or anatomical neck of the humerus. Despite its clinical relevance, the mechanical failure behaviour of long bones is still not sufficiently understood. \textcolor{changes}{To the best of our knowledge, currently, no fully verified and validated numerical methods exist for predicting fractures in the human humerus, which could assist surgeons in medical decisions.}\\
Several previous studies demonstrated the potential of personalized finite element analyses based on quantitative computed tomography (QCT) to predict the strength of human bone. In those studies, the bone is considered as \textcolor{changes}{an elastic isotropic material} with heterogeneous material properties. Commonly, the Young's modulus is a function of the ash density, which can be obtained from QCT scans~\cite{knowles2016quantitative}. For example, the mechanical response of femurs was analysed and compared with in-vitro experiments for intact femurs~\cite{yosibash2007ct,yosibash2007reliable}, femurs with implants~\cite{katz2018patient} and femurs with metastatic tumors~\cite{trabelsi2014patient}. Predictions on the failure of bone are usually based on uncoupled fracture criteria or indicators, which use the computed strains to determine the yield force and the location at which failure starts. Commonly used failure criteria are based on von Mises stress \cite{dragomir2011robust}, yield stress \cite{keller1994predicting, keyak1994correlations}, and maximum principal strain \cite{yosibash2007reliable, yosibash2010predicting, schileo2007subject} and have been used extensively to predict the onset of failure in human femurs. \\
\textcolor{changes}{In-vitro experiments and validated finite element analyses of the proximal humeri are sparse}. Destructive experiments on humeri were performed in \cite{skedros2016radiographic, fankhauser2003cadaveric}, and \cite{dahan2019finite} presented in-vitro experiments inducing physiological fractures in the anatomical neck of the proximal humerus. QCT based finite element analysis was performed and a comparison of strains in the linear elastic range suggested that the $E(\rho)$ relationship used for the femur can be transferred to the humerus. In a follow up study, Dahan et al.~\cite{dahan2020neck} considered the physiologically more common surgical neck fractures. A comparison of Digital Image Correlation (DIC) recorded strains with FE analysis showed a moderate agreement on the humeral neck. Yield load predictions based on cortex failure laws could not reproduce the experimental observations indicating that the fracture initiates inside the humeral head. Thus, advanced failure criteria are needed to predict the failure of human long bones in locations where the cortex is very thin and the failure is not brittle. In addition to the possibility of using more involved non-linear models such as \cite{panagiotopoulou2021experimental}, who suggest to use an elasto-plastic material with a von Mises yield condition to predict the failure load, numerical methods for the modeling of crack initiation and propagation can directly be applied to bone fracture. In contrast to the criteria-based fracture models, these methods have the advantage of not only predicting the yield load or onset of fracture, but can provide information about the type of failure including the fracture path. \textcolor{changes}{Different numerical approaches have been proposed to study fracture of bones at the macro scale.} Hambli \cite{hambli2013robust} coupled a QCT finite element model with a quasi-brittle damage law to predict hip fracture and compare force-displacement curves with experimental data of ten proximal femurs. A finite element based cohesive zone model (CZM) was employed by \cite{ural2006cohesive} to study human cortical bone. Ali et al.~\cite{ali2014specimen} proposed a person-specific model based on the extended finite element method (XFEM) to predict failure loads and in vitro fracture patterns in femurs. However, the bone strength could not be predicted accurately and the approach suffered from convergence problems. Recently, Gustafsson et al.~\cite{gustafsson2021subject} validated an XFEM method on two femurs and reproduced both fracture patterns and bone strength. \\
Phase-field models \textcolor{changes}{(PFM)} for fracture present a promising alternative to the numerical modeling of fracture in a classical finite element setting by reducing mesh-dependency problems and the need for ad-hoc criteria on crack initiation and propagation. Based on a smoothed or regularized approximation of the discrete crack, the propagation of fractures follows from the solution of a minimization problem. PFMs have been applied to brittle fracture in homogeneous and isotropic materials \cite{amor2009regularized, miehe2010thermodynamically}, to ductile fracture \cite{ambati2015phase, borden2018phase} and a wide range of different materials such as polymers \cite{yin2020fracture}, concrete \cite{nguyen2016initiation},  poro-elastic media \cite{aldakheel2021global} and polycrystals \cite{nguyen2017multi}. \textcolor{changes}{Shen et al.~\cite{shen2019novel}} proposed a \textcolor{changes}{PFM} for long bones based on a spatially varying energy release rate and presented an experimental comparison with anatomical neck fracture based on one humerus. After calibrating the phase-field length-scale parameter and the inhomogeneous energy release rate they were able to successfully predict the fracture initiation and propagation. \textcolor{changes}{However, the validity of this study is limited as only one bone is considered, which is used solely for the calibration of parameters. Moreover, the proposed PFM is based on a non-common degradation function.} \\
In this contribution, we \textcolor{changes}{enhance the PFM} by Shen et al.~\cite{shen2019novel} and present a numerical framework for the simulation of fracture in human humeri as well as its full validation on three in vitro experiments. \textcolor{changes}{The approach adapts an inhomogeneous energy release rate (as in Shen et al.~\cite{shen2019novel})} into a phase-field model for brittle fracture and combines it with an embedded domain approach, the Finite Cell Method (FCM) \cite{duster2008finite}. As \textcolor{changes}{shown} in \cite{nagaraja2019phase, hug20203d}, combination of a phase-field model with the FCM allows for a flexible and efficient framework to predict failure in complex geometries without the need of generating boundary conforming meshes. After calibration of the phase-field length-scale parameter based on \emph{one} \textcolor{changes}{humerus}, it is shown that the proposed model is able to replicate the DIC monitored strains, experimental failure loads and crack patterns for \emph{all} humeri.\\
The paper is structured as follows: In Section~\ref{sec:methods}, the experimental methods are introduced, followed by the phase-field approach and its calibration. In Section~\ref{sec:results}, the numerical results are presented including a comparison of the strains on the bones' surface, the failure loads and fracture paths. Finally, in Section~\ref{sec:discussion} the results and limitations of the study are discussed and concluded in Section~\ref{sec:conclusion}.
\begin{figure}[b!]
	\centering
	\includegraphics[width=0.85\textwidth]{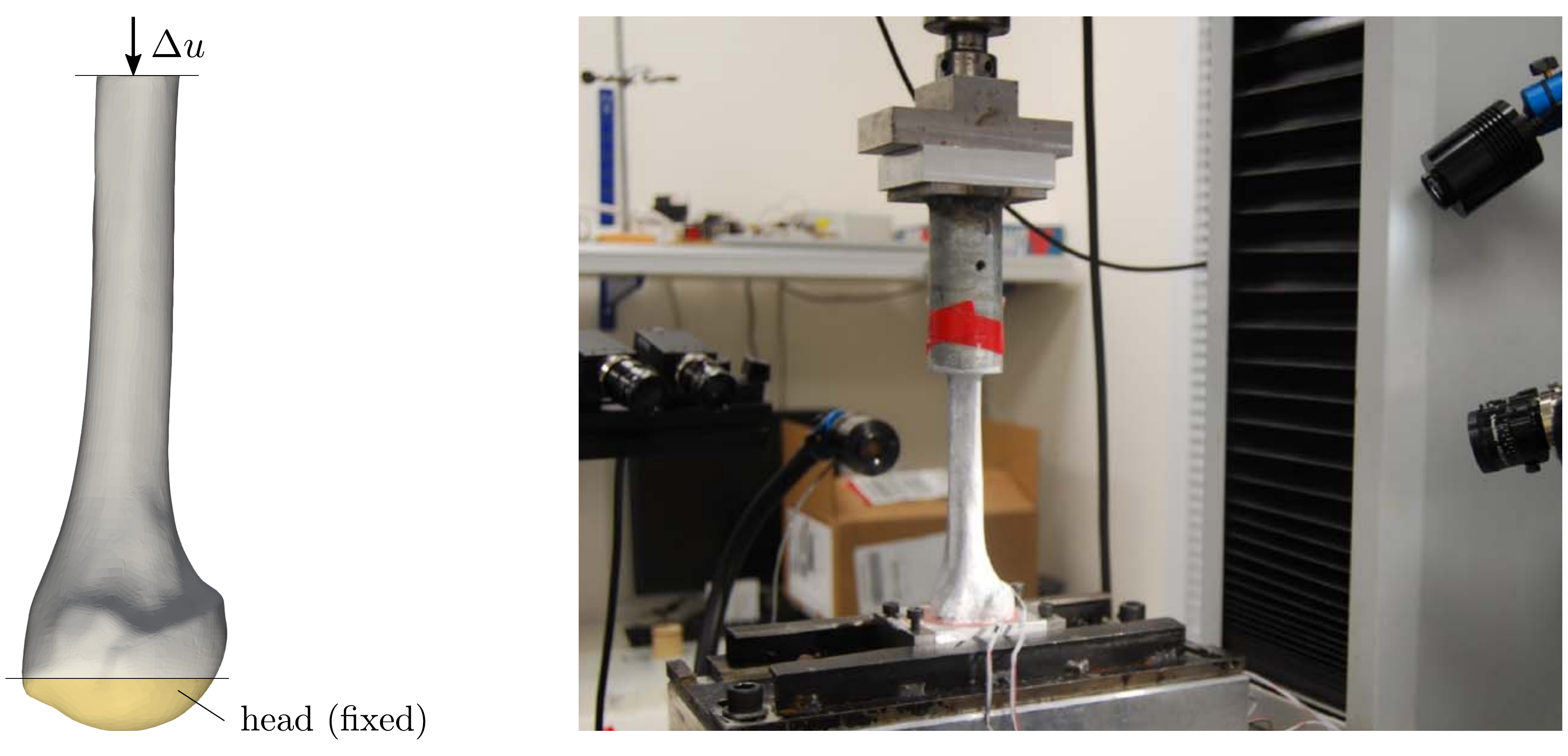}
	\caption{Schematic illustration of \textcolor{changes}{one humerus} in the testing machine with the humeral head embedded in PMMA (left), and photo of the experimental setup with DIC imaging (right) from~\cite{dahan2020neck}.}
	\label{fig:expsetup}
\end{figure}

%% file: sections/methods.tex
\section{Methods} \label{sec:methods}
\subsection{Mechanical Experiments} 
\begin{figure}[b!]
	\centering
	\includegraphics[width=0.75\textwidth]{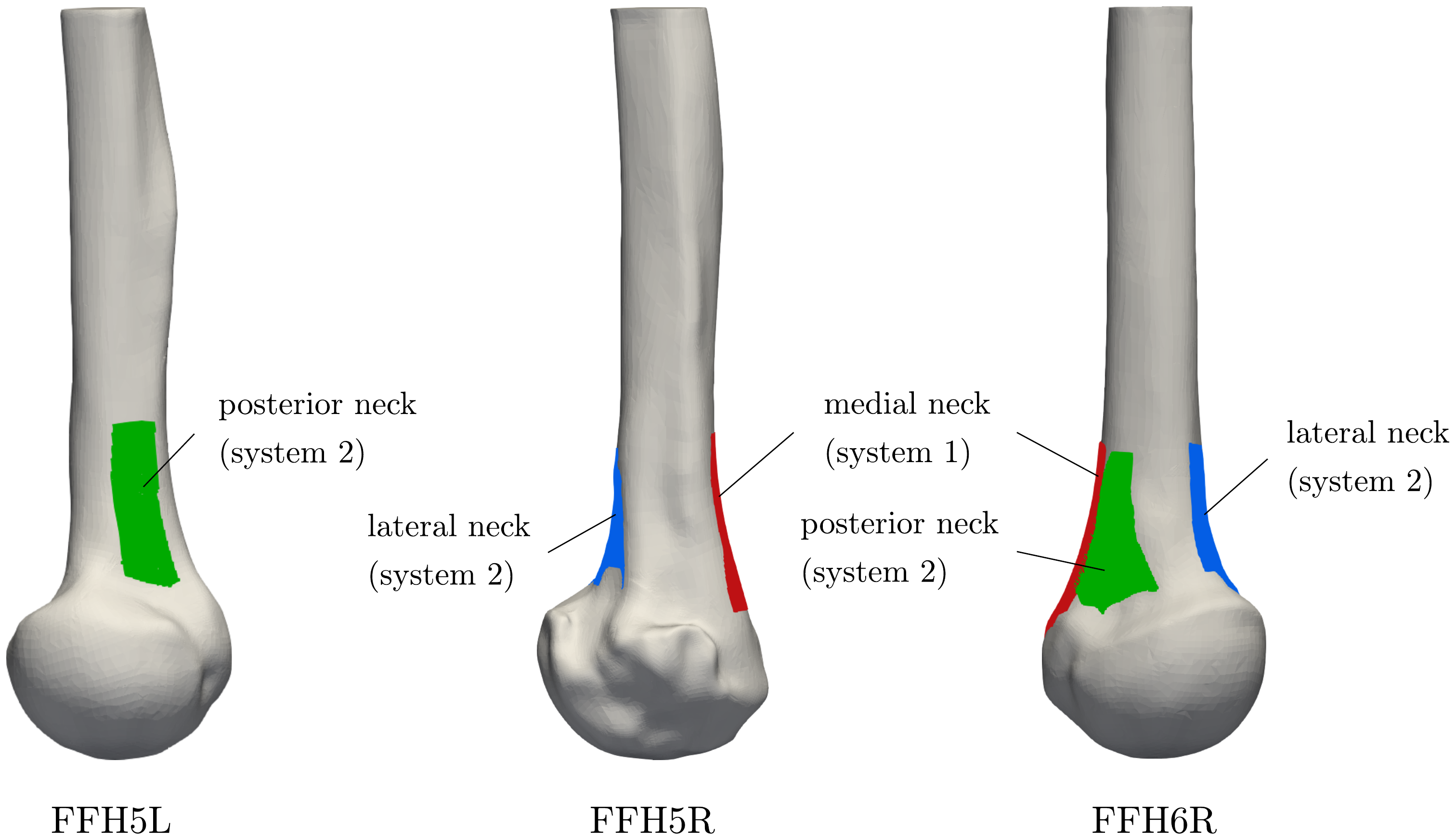}
	\caption{Choice of AOIs for the different humeri after \cite{dahan2020neck}. For FFH5L, the strains on the posterior neck are monitored, while  for FFH5R and FFH6R regions on the lateral and medial neck are chosen.}
	\label{fig:setup+dic}
\end{figure}

The experiments on fresh frozen humeri \textcolor{changes}{are documented in} Dahan \cite{dahan2020neck}. Three proximal humeri, denoted FFH5R, FFH5L and FFH6R, are considered in the following. The humeri were kept frozen at $-80 \degree$ until the day of the experiment, when they were defrosted and the soft tissue was removed. The bones were cut $260$ mm below the top of the humeral head and CT scanned along with five $K_2HPO_4$ calibration solutions in a Brilliance $64$ scanner. The experimental setup is shown in Figure~\ref{fig:expsetup}. To induce fracture at the surgical neck, the bones' head was immersed in PMMA and the humeri were loaded in a testing machine with the proximal part pointing downwards. The humeri were loaded until fracture in a AG-IC, Shimadzu machine (Kyoto, Japan) using a displacement controlled setting. A 6-axis load-cell was used to record the reaction forces. Two DIC systems with two cameras each ($35$ mm lenses and two LED spotlights) were positioned on opposite sites of the humeri to monitor the strains. Areas of interest (AOI) were defined along the anatomical neck (see Figure~\ref{fig:setup+dic}), and the cameras were positioned to create overlap of the camera fields. While bones FFH5R and FFH6R were imaged along the medial and lateral neck, the bone FFH5L was imaged along the posterior neck. DIC images were processed by Vic-3D software (Correlated Solutions Inc.) and smoothed using a spatial Gaussian filter and a time filtering to reduce noise. Operation of the DIC system was checked by comparing the DIC measured strains to 12 strain gauges which resulted in differences less than $10\%$. For more details on the experimental setup and specific parameters see \cite{dahan2020neck}.

\subsection{The phase-field model}
\textcolor{changes}{We follow the numerical framework presented in \cite{nagaraja2019phase, hug20203d}, and combine a phase-field model for fracture with the FCM~\cite{parvizian2007finite}, which allows for a flexible representation of the complex bone geometry. Based on the formulation for fracture in human long bones by Shen~et~al.~\cite{shen2019novel}}, a spatially varying critical energy release rate is introduced to \textcolor{changes}{represent} the heterogeneous failure behavior in the humerus. 
\subsubsection{Governing Equations}
\begin{figure}[b!]
	\centering
	\includegraphics[width=0.73\textwidth]{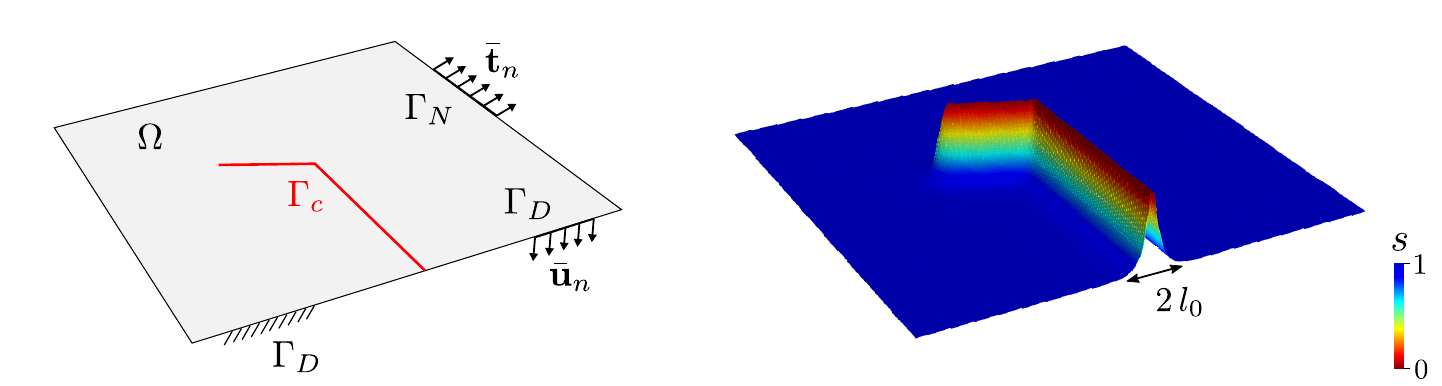}
	\caption{Domain with a discrete crack (left) and its phase-field representation using the continuous phase-field parameter $s$ (right) \cite{hug20203d}.}
	\label{fig:phasefield_setup}
\end{figure}
The phase-field approach to fracture is based on a continuous approximation of the discrete crack using a scalar variable $s$, the so-called phase-field parameter. As \textcolor{changes}{illustrated} in Figure \ref{fig:phasefield_setup}, the phase-field represents the crack as a smooth transition from damaged material ($s = 0$) to intact material ($s=1$) over a regularization width $l_0$. Following the variational formulation by Francfort~et~al.~\cite{francfort1998revisiting} and its regularization by Bourdin~et~al.~\cite{bourdin2000numerical, bourdin2008variational}, crack propagation is formulated as a minimization problem of the functional 
\begin{align}
\label{eq:bourdin_reg_functional}
E_{l_0}({\textbf u}, s)\,=\, \int_{\Omega}\,g(s)\,\Psi(\boldsymbol{\varepsilon}) \textrm{d}{\textbf x}\,+\,\frac{G_{c}}{c_w}\int_{\Omega} \left( \frac{1}{2\,l_0} w(s) \,+\, 2\,l_0 |\nabla s|^2\right) \textrm{d}{\textbf x}\,.
\end{align}
Here, $G_c$ is the critical energy release rate and $g(s)$ is the degradation function which models the loss of stiffness due to damage. We choose a quadratic degradation function $g(s)=(1-\eta)s^2\,+\,\eta$, where $\eta$ is a small parameter which ensures numerical stability when the material is fully damaged. Following the AT-2 model we set the energy dissipation function $w(s)=1\,-\,s^2$ and the scaling parameter $c_w=1/2$ \cite{tanne2018crack}. To prevent crack propagation in compression, a tension-compression split of the elastic strain energy density is commonly used. \textcolor{changes}{Different approaches have been proposed which split the elastic strain energy density into a positive and a negative part, i.e. $\Psi(\boldsymbol{\varepsilon}) = \Psi^+ + \Psi^-$. In Eq.\ (\ref{eq:bourdin_reg_functional}), only the positive part of the elastic strain energy density is degraded by replacing $g(s)\,\Psi(\boldsymbol{\varepsilon})$ with $g(s)\,\Psi^+ + \Psi^-$}. In the following, we use the volumetric-deviatoric split by Amor~et~al.~\cite{amor2009regularized}. Here,
\begin{equation}
\begin{aligned}
&\Psi^+(\boldsymbol{\varepsilon}) = \frac{1}{2}\kappa_0 \langle \textrm{tr}(\boldsymbol{\varepsilon}) \rangle_+^2 + \mu (\boldsymbol{\varepsilon}^\textrm{dev} : \boldsymbol{\varepsilon}^\textrm{dev}) \\
&\Psi^-(\boldsymbol{\varepsilon}) = \frac{1}{2}\kappa_0 \langle \textrm{tr}(\boldsymbol{\varepsilon}) \rangle_-^2
  \end{aligned}
\end{equation}
with the bulk modulus $\kappa_0=\lambda + \frac{2\,\mu}{3}$, the deviatoric strain $\boldsymbol{\varepsilon}^\textrm{dev} = \boldsymbol{\varepsilon} - \frac{1}{3} \textrm{tr}(\boldsymbol{\varepsilon}\,) \boldsymbol{\mathds{I}}$ and the Macaulay brackets $\langle  x \rangle_+ =\frac{1}{2}(x \pm |x|) $. The minimizer of Eq.\ (\ref{eq:bourdin_reg_functional}) is the solution of the associated set of Euler-Lagrange equations. Thus, the strong form of the resulting coupled system of equations reads
\begin{subequations}
	\label{eq:governeq}
	\begin{align}
	\textrm{div}(\bm{\sigma})\,+\,\rho\,{\bm b}&=0,\qquad  \textrm{where } \bm{\sigma} = g(s)\,\dfrac{\partial \Psi^+(\boldsymbol  \varepsilon)}{\partial \boldsymbol \varepsilon} + \frac{\partial \Psi^-( \boldsymbol \varepsilon)}{\partial \boldsymbol\varepsilon} \\[0.1cm]
	-4\,l_0^2\,\Delta s\,+\,(\,\frac{4\,l_0}{G_c}\,(1-\eta)\,\mathcal{H}+\,1)\,s\,&=\,1
	\end{align}
\end{subequations}
and is subject to the boundary conditions
\begin{align}
\label{eq:bc}
{\textbf u} &= \bar{\textbf{u}}_n\quad &&\textrm{on }\Gamma_{D}, \\
{\boldsymbol{\sigma}}\cdot {\textbf{n}} &= \bar{\textbf{t}}_n &&\textrm{on }\Gamma_N, \\
\nabla d \cdot {\textbf{n}} &= 0 &&\textrm{on }\Gamma_{D} \cup \Gamma_N. 
\end{align}
Here, $\mathcal{H}$ is the history variable which replaces the positive part of the elastic strain density~\textcolor{changes}{$ \Psi^+$} in the phase-field equation. Introduced by Miehe~et~al.~\cite{miehe2010thermodynamically} it is defined as
\begin{align}
\mathcal{H}(\textbf{x},t) \coloneqq \underset{t \in [0,T]}{\textrm{max}}\Psi^+(\boldsymbol \varepsilon(\textbf{x},t))\,,
\end{align} 
and ensures irreversibility of the phase-field.
\subsubsection{Numerical Solution}
\begin{figure}[b]
	\centering
	\includegraphics[width=0.75\textwidth]{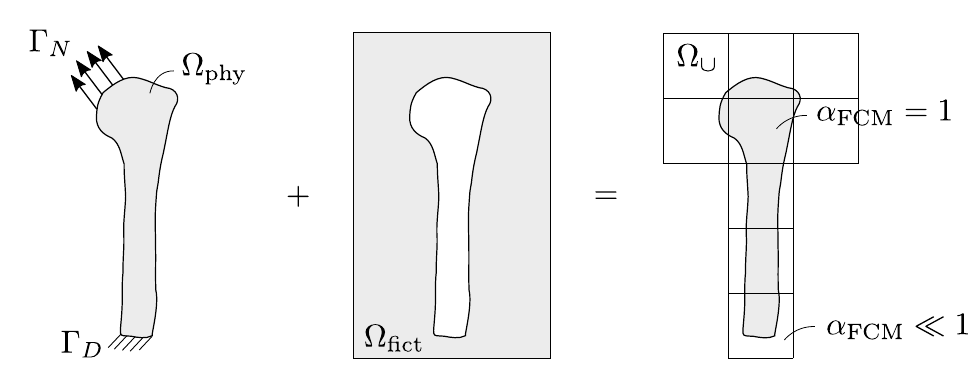}
	\caption{Concept of the FCM following \cite{parvizian2007finite}. Here, cells of the \textcolor{changes}{Cartesian} mesh which do not intersect the domain are not considered as they do not contribute to the stiffness matrix.}
	\label{fig:fcm}
\end{figure}
For the discretization of the coupled system (\ref{eq:governeq}a-b) we use a high-order embedded domain approach, the FCM~\cite{parvizian2007finite,duster2008finite}. By combining the $p$-version of the Finite Element Method (FEM) with an embedded domain approach, the FCM can benefit from high convergence rates while avoiding potentially tedious meshing in the case of complex geometries. As shown in Figure \ref{fig:fcm}, the FCM does not explicitly resolve the original geometry $\Omega_\textrm{phy}$, but embeds it in a fictitious domain $\Omega_\textrm{fict}$ of simple shape. The resulting computational domain $\Omega_\cup = \Omega_\textrm{phy} \cup \Omega_\textrm{fict}$ can be meshed using a structured grid. To recover the original geometry during integration, an indicator function $\alpha(\boldsymbol{x})$ is defined as
\begin{equation}
\alpha({\boldsymbol x})\,=\,\begin{cases} 1\,, \quad \forall {\boldsymbol x} \in \Omega_{phy}, \\ \alpha_{\textrm{FCM}} \ll\, 1\,, \quad \forall {\boldsymbol x} \in \textcolor{changes}{\Omega_{fict}}, \end{cases}
\end{equation} 
which penalizes the contributions of the fictitious domain. Here, $\alpha_{\textrm{FCM}}$ is a very small numerical parameter which avoids ill-conditioning of the stiffness matrix. Special care must be taken when integrating the weak form. The penalization with $\alpha(\boldsymbol{x})$ introduces a discontinuity in the cells cut by the domain boundary and standard numerical integration techniques like Gauss-Legendre quadrature fail to provide accurate results \cite{duster2008finite}. \textcolor{changes}{To overcome this difficulty, we partition the cells into a fine grid of sub-cells, and apply the Gauss-Legendre quadrature for each sub-cell \cite{yang2012efficient}. As geometry and material parameters are also defined voxel-wise this approach is a natural choice for image-based analysis. In contrast to conventional voxel-FEM, the FCM utilizes higher-order shape functions and sub-cell integration techniques and thus is able to use much coarser meshes while providing results of similar accuracy \cite{yang2012efficient}}. Since the boundary of the physical domain does not coincide with the faces of the elements, a penalty method is used to apply Dirichlet boundary conditions \cite{zander2012finite}. Moreover, special care must be taken when computing the reaction forces on the embedded surfaces \cite{d2022accurate}. For the formulation of the weak form the reader is referred to \cite{nagaraja2019phase,hug20203d}. In classic FCM, integrated Legendre polynomials are used as basis functions for the Finite Element test and trial spaces. In this contribution, we use hierarchical B-splines \cite{schillinger2012small}. As proposed e.g.\ in \cite{miehe2010thermodynamically}, the coupled quasi-static problem is solved in an alternating manner using a staggered scheme. In each displacement step, the phase-field and the elastic equation are solved subsequently until the residual drops below a given tolerance $\varepsilon_\textrm{stag}$, i.e. the iterations are terminated after staggered step $i$ if
\begin{equation}
	s_{\textrm{tol}, i} < \varepsilon_\textrm{stag}\,, \quad \textrm{ where } s_{\textrm{tol}, i} = \textrm{max}(\mathcal{R}_{\boldsymbol u}, \mathcal{R}_{s} )\,,
\end{equation}
where $\mathcal{R}_{\boldsymbol u}$ is the residual of the elastic problem and $\mathcal{R}_{s}$ is the residual of the phase-field problem.

\subsubsection{Simulation Setup}
\begin{figure}[b!]
	\begin{floatrow}
		\ffigbox[8.5cm]{%
			\includegraphics[width=0.45\textwidth]{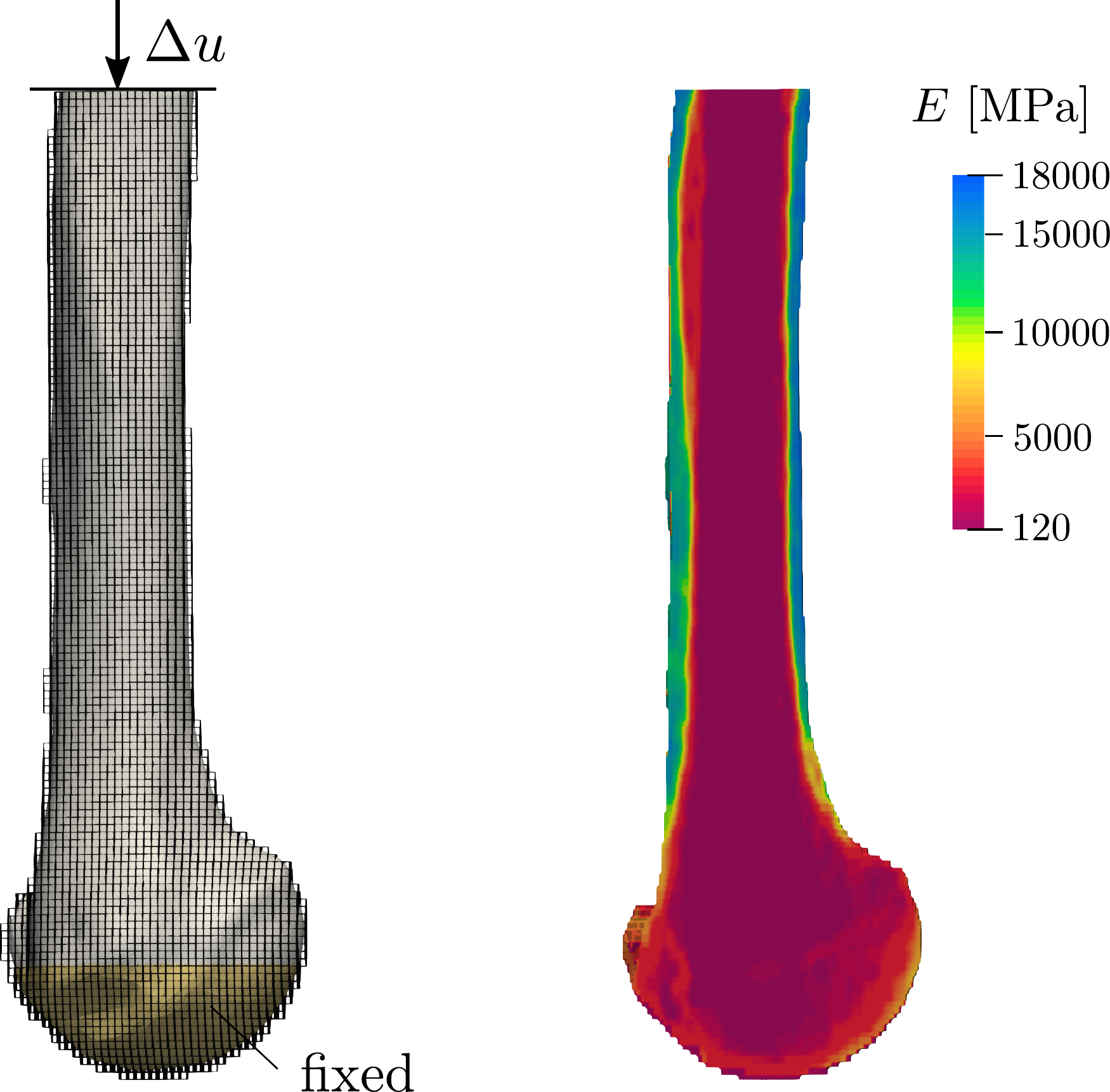}
			\vspace{0.5cm}
		}{%
			\caption{Geometry with boundary conditions and FCM mesh (left), and Young's modulus distribution in the humeri (right).}%
			\label{fig:simsetup}
		}\hfill
		\capbtabbox[7cm]{%
			\setlength{\tabcolsep}{5pt}
			\renewcommand*{\arraystretch}{1.4}
			\begin{tabular}{cc} \hline
				$h$ [mm] & $1.25$ \\
				$p$ [-]& $3$ \\
				$\alpha_{\textrm{FCM}}$ [-]& $1.0\,\cdot 10^{-6}$ \\
				$u_\textrm{large}$ [mm]   &$0.04$  \\
				$u_\textrm{med}$ [mm]  &$0.002$  \\
				$u_\textrm{small}$ [mm]  &$0.001$  \\
				$\varepsilon$ [-]& $1.0\,\cdot 10^{-5}$ \\
				$n_\textrm{stag}$ [-] & 25 \\
				$\nu$ [-] & $0.3$ \\
				$E_0$ [MPa] & 20000 \\
				$G_{c_0}$ [N/mm] & 7 \\ 
				\hline
				\vspace{0.05em}
			\end{tabular}
		}{%
			\caption{Simulation parameters for the three humeri.}
			\label{tab:material}
		}
	\end{floatrow}
\end{figure}
Based on QCT scans of the humeri segmentations are generated following~\cite{dahan2020neck}. The computational domain is defined as the bounding box of the segmentation and discretized using a structured grid. Elements completely inside the fictitious domain are excluded from the model. The boundary conditions are applied \textcolor{changes}{in a weak sense} based on a surface triangulation of the humerus, which is intersected with the computational mesh for accurate integration~\cite{elhaddad2018multi}. As shown in Figure \ref{fig:simsetup}, the humeral head, \textcolor{changes}{embedded in PMMA}, is fixed ($u_x = u_y = u_z = 0$), while a displacement $\Delta u_z$ is applied on the distal face. The young's modulus is computed voxel-wise based on the Hounsfiled units (HU) as described in \cite{yosibash2014predicting,katz2019scanner}. First, HU values are converted to equivalent mineral density $\rho_{K_2HPO_4}$ and then to ash density $\rho_\textrm{ash}$ following~\cite{goodsitt1992conversion,schileo2008accurate} as
\begin{align}
\rho_\textrm{ash} = 0.877 \times 1.15 \times \rho_{K_2HPO_4} + 0.08\,.
\end{align}
Thereafter, Young's modulus values are computed as a function of $\rho_\textrm{ash}$ based on the relations proposed by \cite{keyak1994correlations, keller1994predicting} as
\begin{align}
E_\textrm{cor}(\rho_\textrm{ash}) &= 10200\,\rho_\textrm{ash}^{2.01}\:\textrm[MPa]\,, \qquad &\rho_\textrm{ash} > 0.486\,,&&\\
E_\textrm{tra}(\rho_\textrm{ash}) &= 2398\:\textrm[MPa]\,, \qquad & 0.3 < \rho_\textrm{ash} \leq 0.486\,,&&\\
E_\textrm{tra}(\rho_\textrm{ash}) &= 33900\,\rho_\textrm{ash}^{2.2}\:\textrm[MPa]\,, \qquad &\rho_\textrm{ash} \leq 0.3\,.&&
\end{align}
\textcolor{changes}{Following~\cite{shen2019novel}, the critical energy release rate $G_c$ is assumed to depend on the bone density following a power-law relation based on the Young's modulus
\begin{equation}
G_c\left(\rho_\textrm{ash} \right) = G_{c_0}\,\left(\frac{E(\rho_\textrm{ash})}{E_0}\right)^\beta\,,
\end{equation}}
where $E_0$ and $\beta$ are the base Young's modulus and the power-law exponent. For $\beta=1$, a linear dependence between Young's modulus and the energy release rate is obtained, while for $\beta<1$ less correlation between the two material properties is assumed. Shen et al.~\cite{shen2019novel} investigated the influence of $\beta$ and calibrated \textcolor{changes}{it} to $\beta=0.8$, which we adopt for our computations. \textcolor{changes}{The base values are set to $G_{c_0} = 7$ N/mm~\cite{yeni1998influence} and $E_0 = 20000$ MPa based on different studies on the Young's modulus of human cortical bone \cite{rho1993youngs, hoffmeister2000anisotropy, zysse1999elastic}}. Further simulation parameters are listed in Table~\ref{tab:material}. The computational mesh consists of $77518$ (FFH5R), $81318$ (FFH6R) and $75052$ (FFH5L) finite cells with an edge length of $1.1$ mm, which allows to resolve phase-field length scales $l_0 \geq 1.1$ mm. The length parameter $l_0$ is calibrated in the next Section~\ref{sec:calib}. If not stated otherwise, we use a polynomial degree $p=3$ and penalize contributions of the fictitious domain with $\alpha_{\textrm{FCM}} = 1.0\,\cdot 10^{-6}$. To accurately resolve the crack initiation and propagation, we adapt the size of the displacement steps throughout the simulation. Starting with larger steps $u_\textrm{large} = 0.04$ mm, we first decrease the step size to $u_\textrm{med} = 0.002$ mm, and then to $u_\textrm{small} = 0.001$ mm closer to crack initiation. The tolerance for the staggered solution scheme is set to $\varepsilon_\textrm{stag} = 1.0\,\cdot 10^{-5}$. In each displacement step, a maximum number of $n_\textrm{stag} = 25$ staggered steps is performed. The poisson ratio is set to $0.3$.

\subsection{Calibration} \label{sec:calib}
\begin{figure}[b!]
	\begin{minipage}[b]{0.65\textwidth}
		%\centering
		\includegraphics[width=\textwidth]{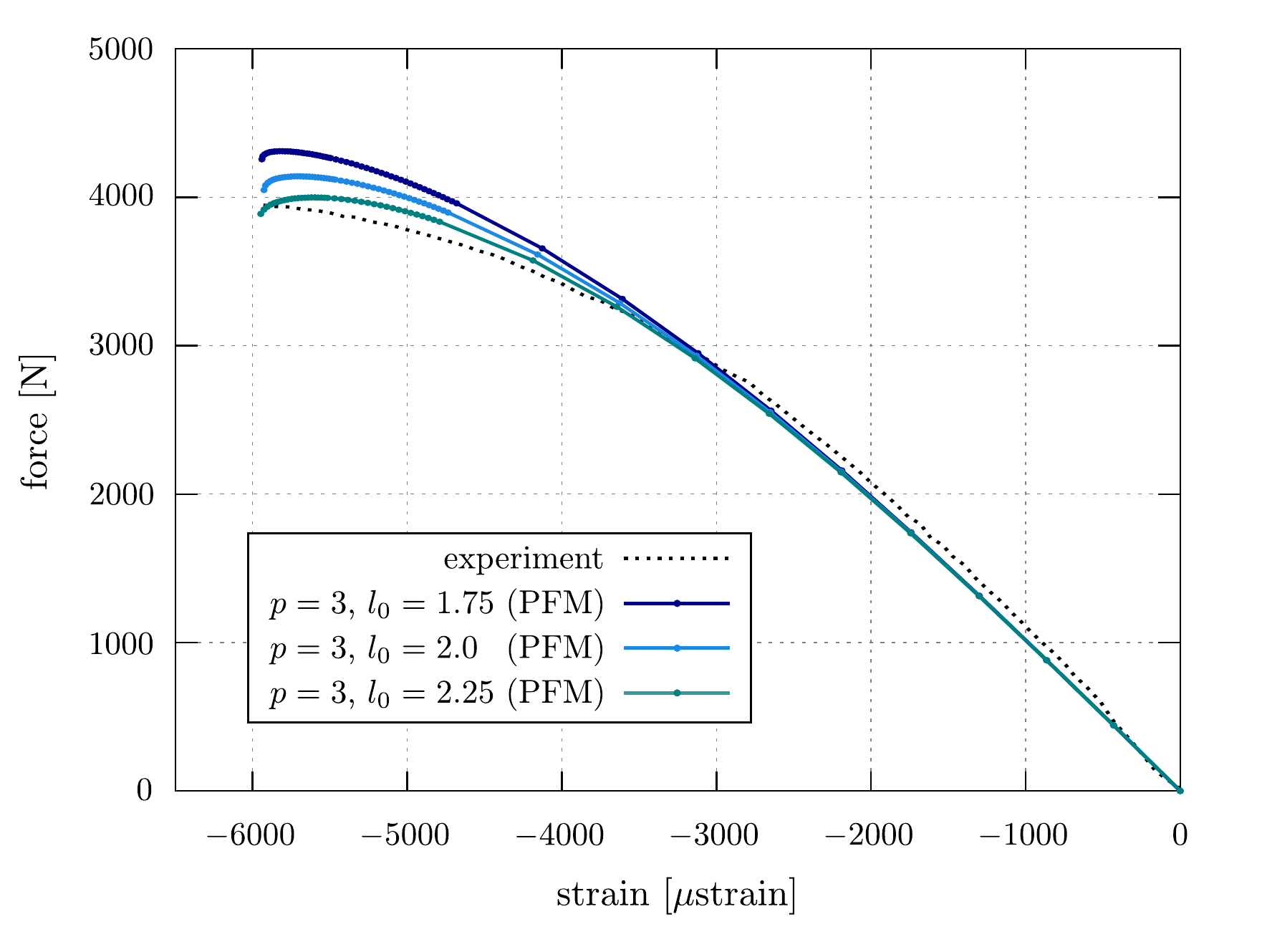}
	\end{minipage}%
	\begin{minipage}[b]{0.35\textwidth}
		\small
		\hfill
		\centering
		\setlength{\tabcolsep}{5pt}
		\renewcommand*{\arraystretch}{1.4}
		\begin{tabular}{p{1.85cm}p{1.5cm}p{1.5cm}}
			&   failure load [N]   & rel.$\qquad$   error [$\%$] \\ \midrule
			experiment & $3948$ & - \\
			$l_0=1.75$ & $4310$ & $9.2$ \\
			$l_0=2.0$ & $4141$ & $4.8$ \\ 
			$l_0=2.25$ & $3999$ & $1.3$\\
		\end{tabular}
		\vspace{4.5cm}
	\end{minipage}
	\caption{Calibration of the length-scale parameter $l_0$ based on the humeri FFH5R. Force-strain curves at a point close to fracture location for different choices of $l_0$ in comparison with the experimentally measured curve (left). Comparison of failure loads (right).}
	\label{fig:l0_comp}
\end{figure}
\begin{figure}[t!]
	\centering
	\includegraphics[width=0.98\textwidth]{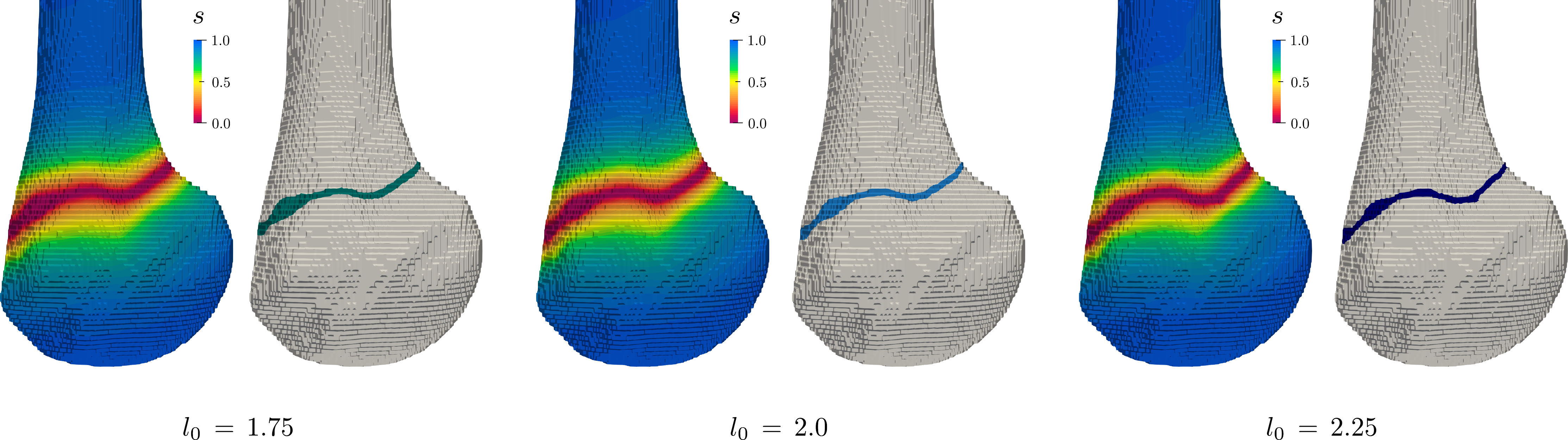}
	\caption{Comparison of crack patterns for different length-scale parameter $l_0$ for the humeri FFH5R.}
	\label{fig:l0_comp_visu}
\end{figure}

An integral part of phase-field models is the length-scale parameter $l_0$. Originally introduced as a purely numerical parameter, which recovers the discrete nature of fracture in the limit $l_0 \rightarrow 0$, several contributions interpret $l_0$ as a material parameter~\cite{del2013variational, freddi2010regularized, tanne2018crack}. Instead of choosing $l_0$ as small as possible the length-scale parameter is calibrated to match the critical strength of the material. Following this approach we calibrate $l_0$ based on the humeri FFH5R to match the force-strain curve recorded in the experiment. The experimental curve was obtained by plotting the measured reaction force over the principal compression strain $\varepsilon_3$ obtained from the DIC results. The strains were extracted at a location close to where the fracture initiates~\cite{dahan2020neck}. The corresponding computed force-strain curves are shown in Figure~\ref{fig:l0_comp} for different values of $l_0$. A larger length-scale parameter $l_0$ results in a lower material strength and consequently a lower failure load. At the same time the absolute value of the principal compression strain at failure decreases. Due to the increased region of damage more material is degraded as soon as the fracture initiates, which results in an earlier and more pronounced deviation from the linear elastic slope. For all chosen length-scale parameters, the shape of the experimental force-strain curve can be captured. A quantitative comparison of the relative errors in the failure loads is presented in Figure~\ref{fig:l0_comp}, right. Here, the failure load is obtained as the maximum force in the force-strain curve. The best match of computed and experimental failure loads is obtained for $l_0=2.25$~mm with a deviation in the failure load of $1.3\%$. For $l_0=2.25$~mm the error is slightly higher with $4.8\%$ and still below $10\%$ if we choose a length-scale of $l_0=1.75$~mm. In Figure~\ref{fig:l0_comp_visu}, a visual comparison of the obtained crack patterns is presented. For each length-scale parameter, the computed phase-field in the final displacement step is shown. For a clearer visualization of the resulting crack path an iso-volume of the phase-field following $0 \leq s \leq 0.03$ is extracted which corresponds to the fully broken region. Varying $l_0$ in the considered range of values has no decisive impact on the failure pattern and the resulting crack path. Following those observations, we will set the length-scale parameter $l_0$ to $2.25$ mm for all humeri.

%% file: sections/results.tex
\section{Results} \label{sec:results}
As a first step of the analysis, the principal strains in the linear elastic range are validated in Section~\ref{sec:strains}. FCM computed strains are compared with the DIC measured values in a linear regression analysis and a Finite Element reference solution by Dahan~\cite{dahan2020neck}. The R$^2$ and root mean square error (RMSE) are evaluated for all humeri. In Section~\ref{sec:loads}, the experimental and numerical force-strain curves are compared for all humeri. Numerical strains are obtained by averaging over a sphere with radius $0.5$ mm at the location where the DIC strains are measured. The shape of the load-strain curves is compared qualitatively and the average error in failure loads is computed. The computed crack patterns are analysed in Section~\ref{sec:patterns}. Here, photos of the fractured humeri are used for a qualitative comparison of the PFM-FCM predicted fracture patterns, and initiation of the neck fractures is investigated.

\begin{figure}[b!]
	\centering
	\includegraphics[width=1.0\textwidth]{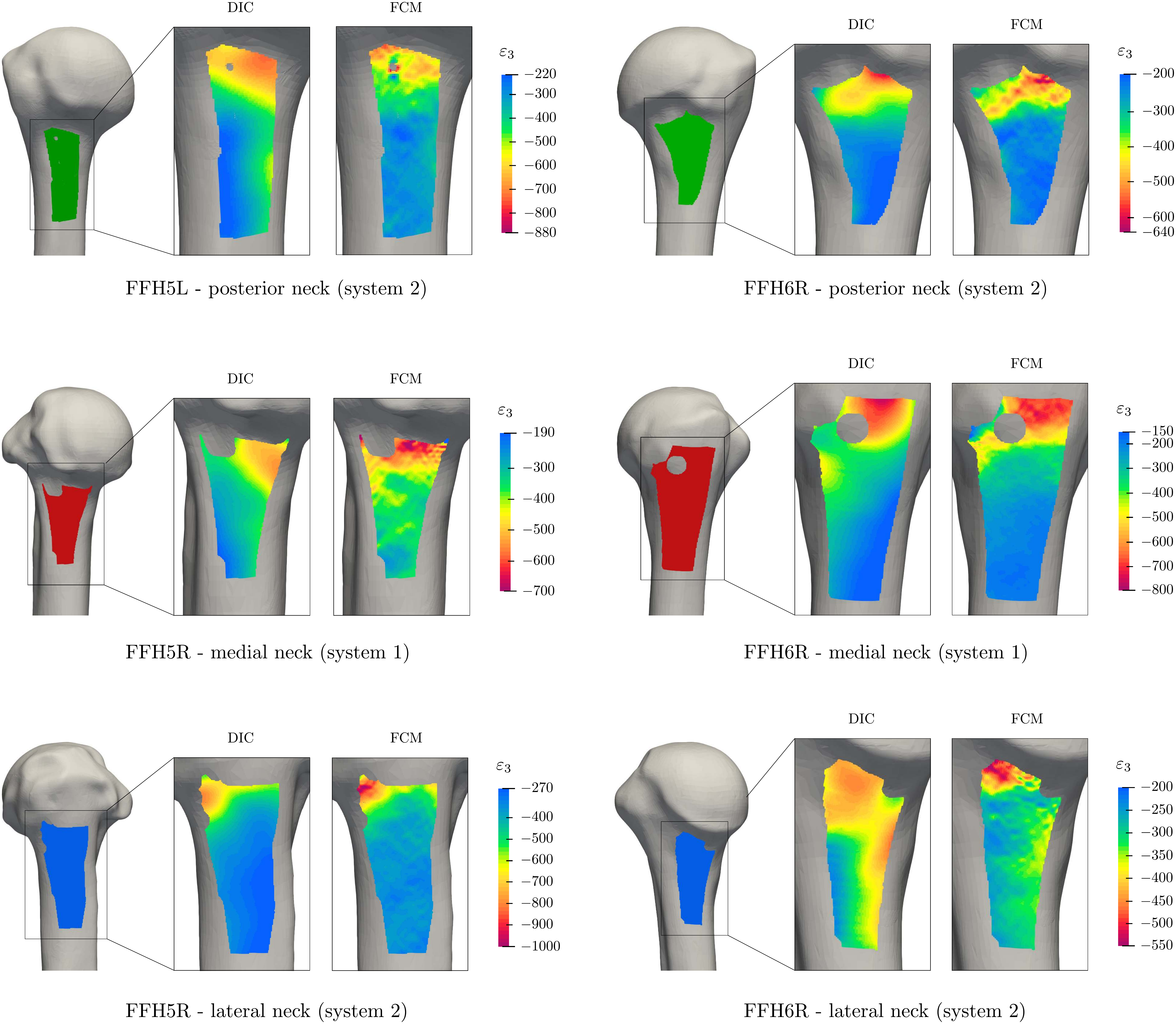}
	\caption{Qualitative comparison of DIC recorded and FCM computed third principal strains for the humeri FFH5L, FFH5R and FFH6R.}
	\label{fig:DIC_humeri_visu}
\end{figure}

\subsection{Strain Validation} \label{sec:strains}
\begin{figure}[h!]
	\centering
	\includegraphics[width=0.98\textwidth]{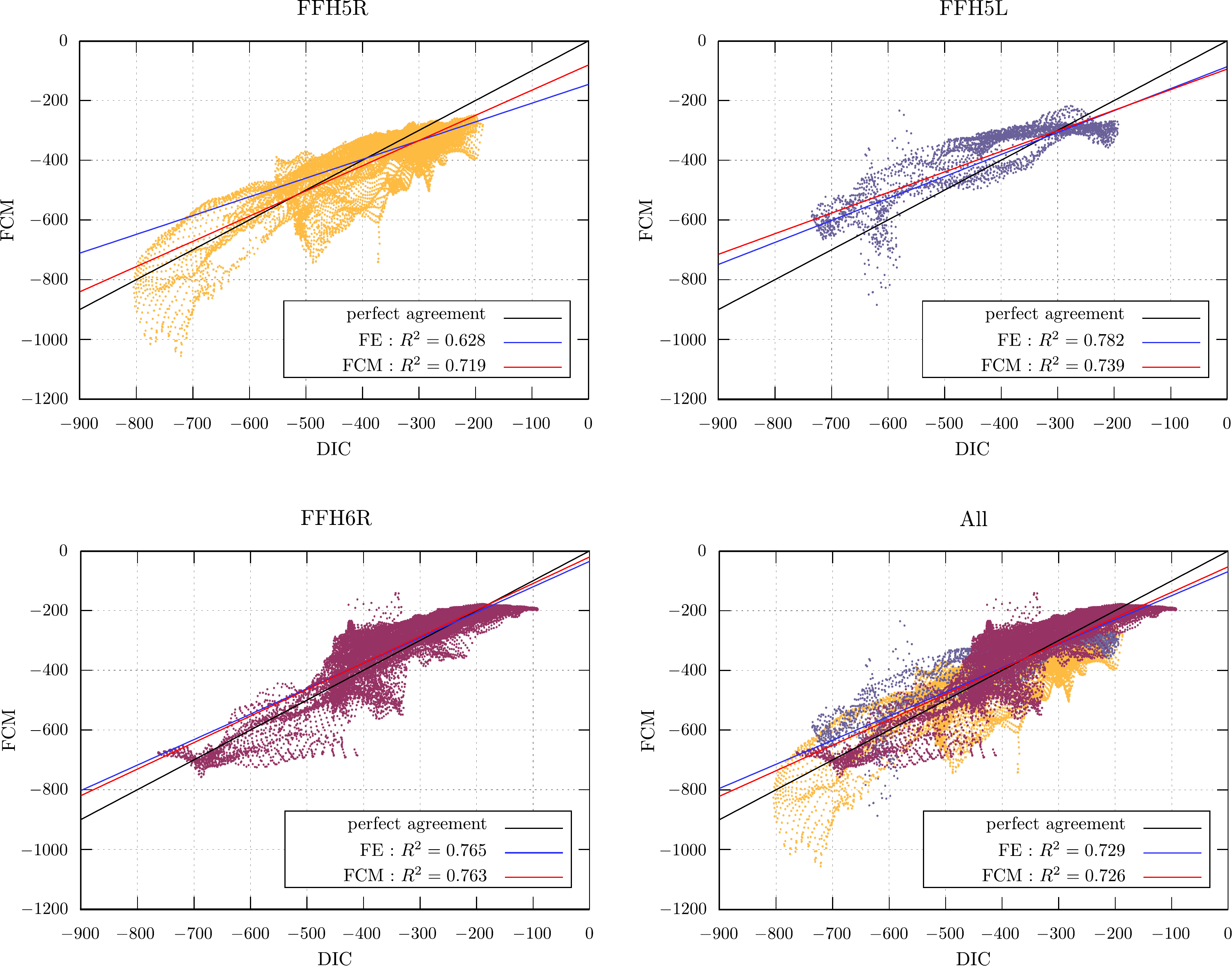}
	\caption{Comparison of the FCM computed principal compression strains with the DIC measurements in a linear regression analysis. Trend lines of the FE reference solution by Dahan~\cite{dahan2020neck} are included in all plots and R$^2$ values are given.}
	\label{fig:DICcomp_scatter}
\end{figure}

For a validation in the linear elastic scheme, strains on the humeral neck are recorded with DIC for all three humeri. In all experiments, a load of $600$ N is applied and principal compression strains $\varepsilon_3$ are obtained. A qualitative comparison of the monitored and computed compression strains is presented in Figure \ref{fig:DIC_humeri_visu}. Clearly, the overall strain distribution can be captured. However, the analysis tends to underestimate $\varepsilon_3$ in regions with minimum compression strains. It should be noted, that in contrast to the DIC values the computed principal strains are not smoothed. Linear regression plots for all three humeri are presented in Figure \ref{fig:DICcomp_scatter}, and computed coefficients and evaluation quantities are listed in Table \ref{tab:DICresults}. The analysis confirms, that the principal strains can be reproduced with reasonable accuracy and are of similar quality for all three humeri with $R^2_\textrm{FCM}=0.719$ (FFH5R), $R^2_\textrm{FCM}=0.739$ (FFH5L) and $R^2_\textrm{FCM}=0.763$ (FFH6R). Results obtained with the FCM are of comparable accuracy to the FE reference solution with $R^2_\textrm{FCM}=0.726$ and  $R^2_\textrm{FE}=0.729$, respectively. Considering all three humeri, the average RMSE is $65.9$ $\mu$strain and average percentage error is $15.4\%$ (RMSE$_\textrm{FE}= 66$~$\mu$strain, $\% e_{rel,\textrm{FE}} =15\%$).

\begin{table}[h!]
	\vspace{0.5cm}
	\small
	\setlength{\tabcolsep}{7pt}
	\renewcommand*{\arraystretch}{1.4}
	\begin{tabular}{lcccccccccc} \toprule
		& \multicolumn{2}{c}{slope} & \multicolumn{2}{c}{intercept} & \multicolumn{2}{c}{$R^2$} & \multicolumn{2}{c}{RMSE} & \multicolumn{2}{c}{$\% e_{rel}$} \\
		& FE & FCM 	& FE & FCM	& FE & FCM	& FE & FCM & FE & FCM \\ \midrule
		FFH5R & $0.628$ & $0.844$ & $-146.0$ & $-80.8$ & $0.628$ & $0.719$ & $71.0 $ & $60.1$ & $15.7$ & $15.9$\\
		FFH5L & $0.736$ & $0.689$ & $-86.5$\textcolor{white}{0} & $-94.7$ & $0.782$ & $0.739$ & $70.0$ & $58.7$ & $15.0$ & $17.0$  \\ 
		FFH6R & $0.853$ & $0.888$ & $-35.1$\textcolor{white}{0} & $-20.4$ & $0.765$ & $0.763$ & $60.0$ & $58.4$ & $15.4 $ & $14.5$ \\
		\bottomrule
	\end{tabular}
	\vspace{0.1em}
	\caption{Linear regression coefficients for the comparison of DIC with FE (reference) and FCM computed strain fields.}%
	\label{tab:DICresults}
\end{table}
\subsection{Failure Loads} \label{sec:loads}
\begin{figure}[t!]
	\centering
	\includegraphics[width=0.75\textwidth]{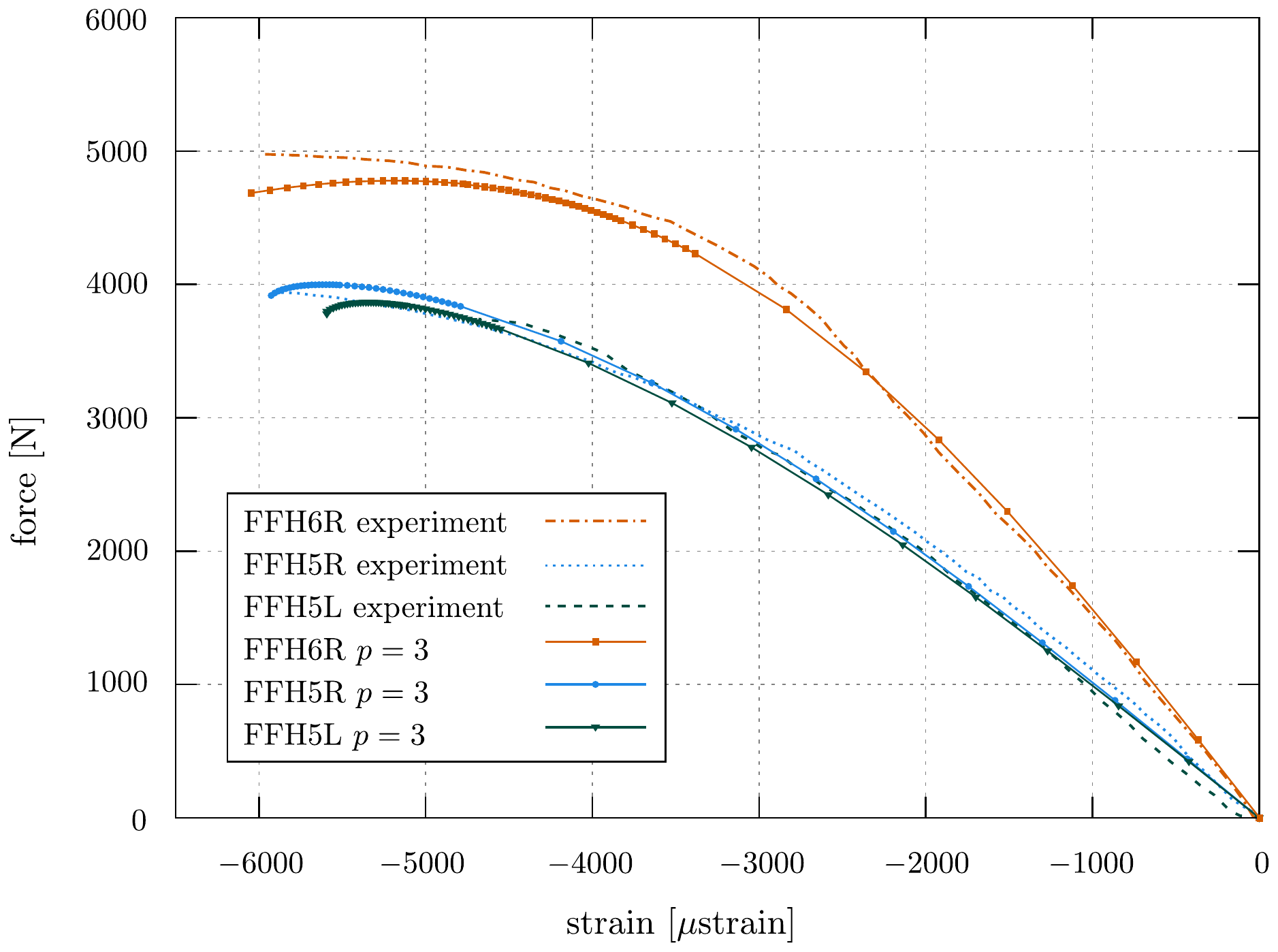}
	\caption{Analysis of the failure behaviour in terms of force-strain curves for the calibrated length-scale $l_0=2.25$ mm. Force-strain curves at a point close to fracture location for the three humeri in comparison with the experimentally measured curve.}
	\label{fig:ForceVsStrain}
\end{figure}

The computed force-strain curves are shown in Figure~\ref{fig:ForceVsStrain} for all three humeri. Curves are obtained from the PFM-FCM simulation by averaging $\varepsilon_3$ over a sphere with radius $0.5$~mm. Clearly, the numerical model is able to capture the fracture behaviour of the three humeri visible in the force-strain curves. In addition to the calibrated humerus FFH5R, but also for the other two humeri the linear elastic range, the onset of fracture and the failure point can be reproduced. For FFH5L, the compression strain at failure is underestimated, while a very good agreement is obtained for the other two humeri. A quantitative comparison of the failure loads is presented in Table \ref{tab:FailureLoads}. The failure loads can be reproduced with very good accuracy for all three humeri with a relative error of $0.6\:\%$ for FFH5L, \textcolor{changes}{$3.8\:\%$ for FFH6R} and $1.3\:\%$ for the calibrated bone FFH5R. The numerical simulation tends to slightly overestimate the failure load: in the case of FFH5R by $60$~kN, for FFH5L with $42$~kN and FFH6R with $53$~kN.

\begin{table}[h!]
\small
\centering
\setlength{\tabcolsep}{5pt}
\renewcommand*{\arraystretch}{1.4}
\begin{tabular}{lp{1.4cm}p{1.4cm}p{1.5cm}p{1.4cm}}
	\centering
	bone & exp. load [kN]    & PFM load [kN] & abs. $\quad$ error [kN] & rel. $\quad$ error [$\%$] \\ \midrule
	FFH5L & $3736$  & $3757$ & $42$ & $0.6$  \\ 
	\textcolor{changes}{FFH6R} &  $4970$ & $4781$ & $189$ & \textcolor{changes}{$3.8$} \\ 
	FFH5R$\:\:$ & $3948$  & $3999$ & $51$ & $1.3$ \\ 
\end{tabular}
	\vspace{0.1em}
	\caption{Comparison of failure loads for the humeri FFH5L and FFH6R, as well as the bone used for calibration (FFH5R).}%
	\label{tab:FailureLoads}
\end{table}
\begin{figure}[h!]
	\centering
	\vspace{0.5cm}
	\includegraphics[width=0.98\textwidth]{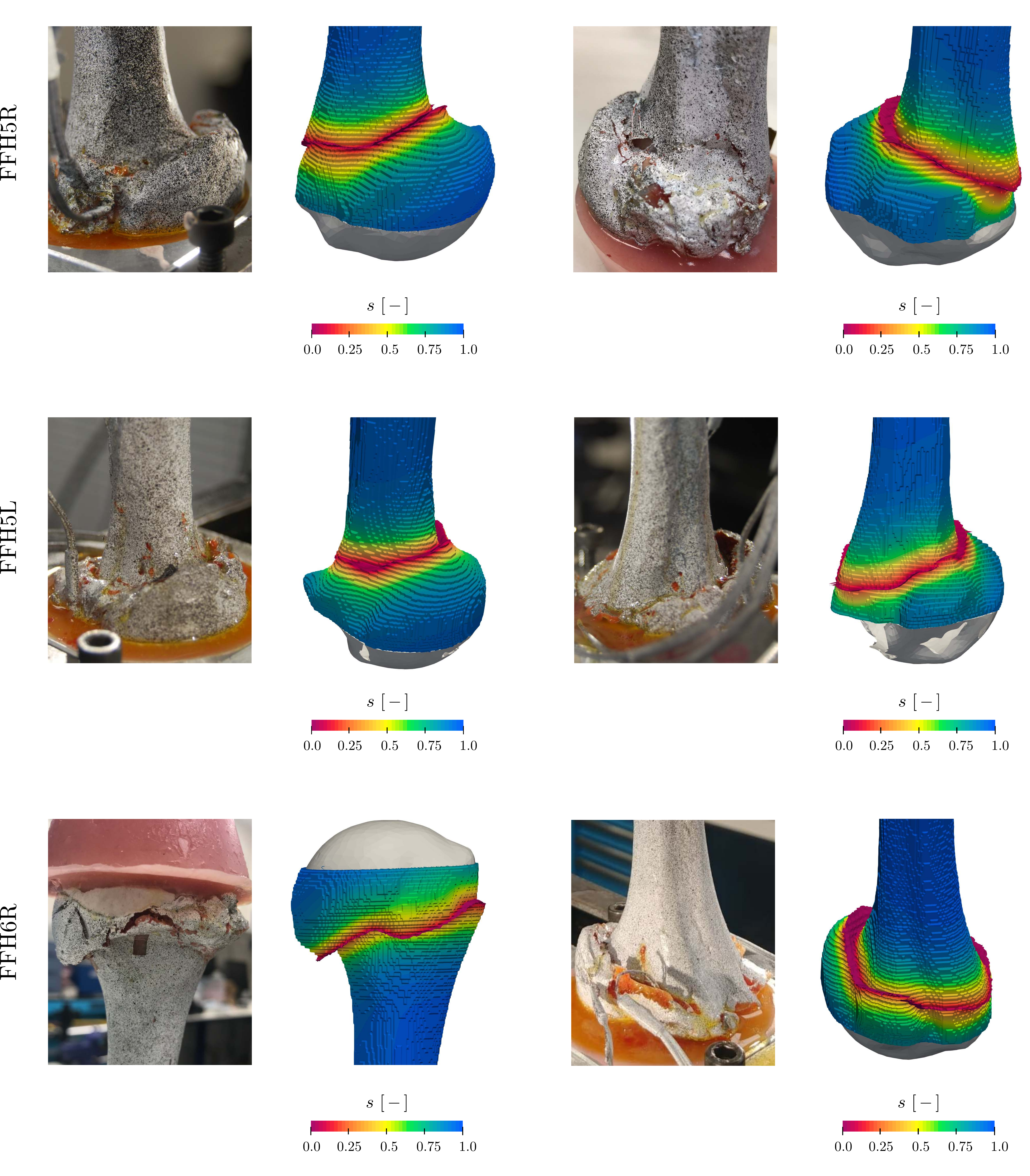}
	\caption{Qualitative comparison of computed crack patterns with experimental results for the three humeri FFH5R, FFH5L and FFH6R. Pictures of the fractured humerus are contrasted with the computed phase-field value. The latter is plotted on an iso-volume of the phase-field $0.035 \leq s \leq 1.0$ which is deformed using the computed displacements magnified by a factor of 8. The grey surface corresponds to the part of the humeral head where the fixed boundary condition is applied.}
	\label{fig:CrackPattern}
	\vspace{0.5cm}
\end{figure}

\begin{figure}[h!]
	\centering
	\includegraphics[width=0.9\textwidth]{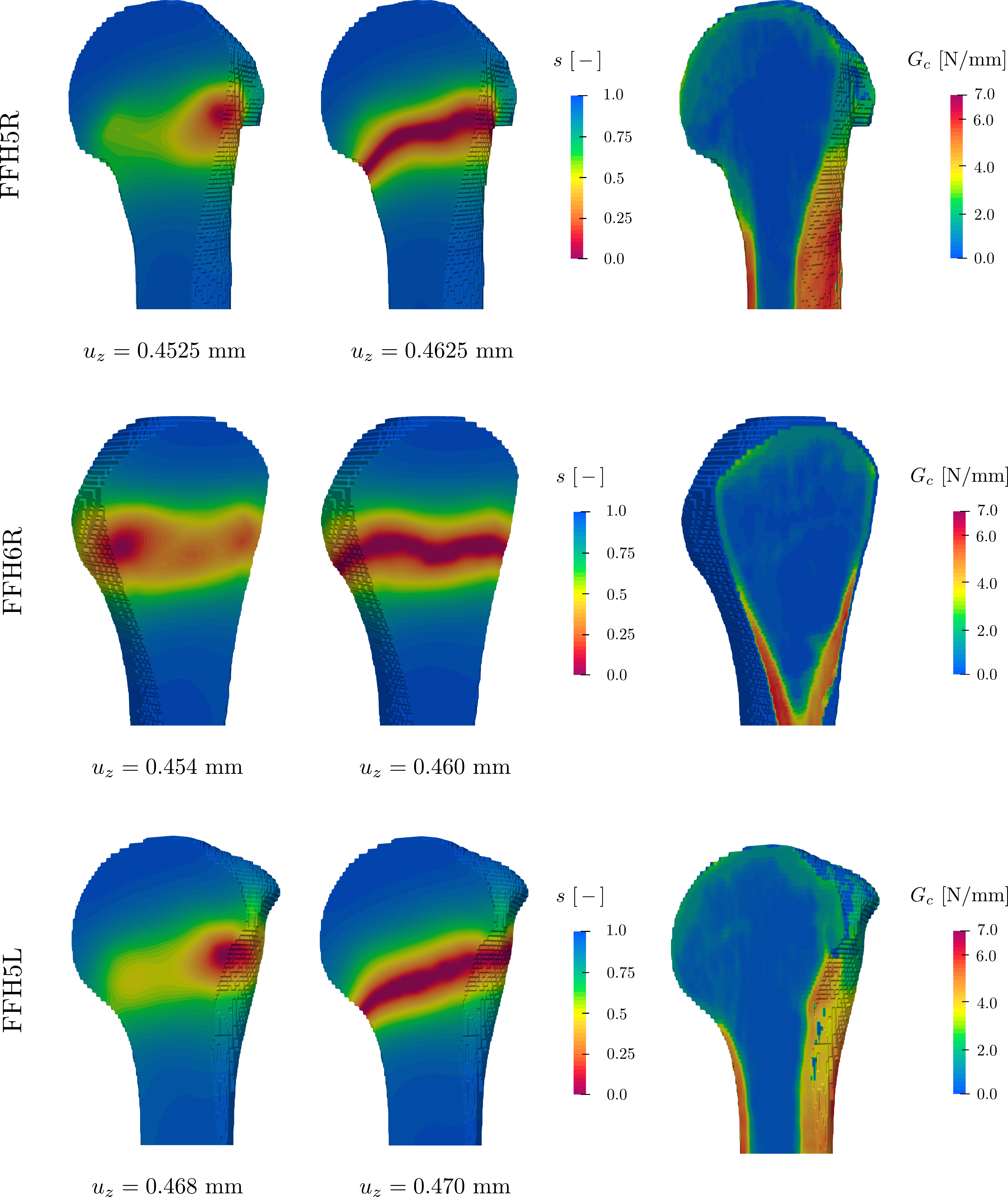}
	\caption{Phase-field parameter $s$ at the point of initiation of the fracture and when the crack is fully developed (left), and critical energy release rate $G_c$ (right) on a cut through the geometry for the three humeri FFH5R, FFH6R and FFH5L.}
	\label{fig:CrackInit}
\end{figure}

\subsection{Crack Pattern} \label{sec:patterns}
To evaluate the computed crack patterns photos of the fractured humeri are compared with the numerical results in Figure \ref{fig:CrackPattern}. For the visualization, an iso-volume of the phase-field $0.035 \leq s \leq 1.0$ is extracted and the geometry is deformed using the computed displacements magnified by a factor of 8. For each humerus, the phase-field parameter is plotted on the geometry. The part of the humeral head that is kept fixed during both the simulation and the experiment is depicted in grey. The qualitative comparison shows, that the fracture at the surgical neck of the humerus can be reproduced for all experiments. Due to the speckle pattern painted on the humeri for the DIC algorithms it is difficult to extract the exact fracture path from the photos. However, the path of the crack observed experimentally clearly shows characteristic features that are captured in the numerical simulation, such as the kink observed for the humerus FFH6R. Similar features can also be found in the humerus FFH5L. In Figure \ref{fig:CrackInit}, the initiation of the fractures as observed with the PFM is visualized. The phase-field is shown at the displacement where the fracture initiates and when the crack is fully developed on a cut through the geometry. Additionally, the right image shows the critical energy release rate on the same geometry. For all humeri, the fracture initiates at the surgical neck of the humerus. In the case of FFH5L and FFH5R, the point of initiation lies on the outer shell of the humeral head. However, damage to the trabecular structures inside the bone can already be observed far away from the location of initiation. For FFH6R, failure is initiated at several points inside the trabecular bone, which upon loading connect to a single fracture extending to the outer shell of cortical bone.

%% file: sections/conclusions.tex
\section{Discussion} \label{sec:discussion}
The validation of the proposed phase-field FCM model is based on three in vitro experiments on surgical neck fractures in human humeri \cite{dahan2020neck}. In a first step, the phase-field FCM model is calibrated based on one of the three experiments. To this end, the length-scale parameter $l_0$ is chosen to match the computed to the experimentally measured force-strain curves of the humeri FFH5R. The best agreement is obtained for $l_0 = 2.25$ mm with a relative error in the failure load of $1.3\:\%$. Varying $l_0$ in the range of $1.75$~mm to $2.25$~mm leads to no visible change in the fracture path and an acceptable error in failure loads below $10\:\%$ is obtained. The results indicate that in the presented setting, the phase-field FCM model is relatively robust with respect to $l_0$.\\
The validation in the linear elastic regime is based on a comparison of the principal compression strains on the humeral neck and shows that the overall experimental strain pattern can be reproduced for all three humeri. Regions with maximum and minimum $\varepsilon_3$ agree particularly well for the posterior and medial neck of FFH6R, while for FFH5L and FFH5R larger differences are visible. These deviations can partly be attributed to the fact that raw numerical values are compared with smoothed DIC strains. A linear regression analysis confirms, that experimental strains can be reproduced with sufficient agreement and that the results for the different humeri FFH5R, FFH5L and FFH6R are of similar accuracy. Considering all three humeri, the correlation of FCM computed strains with experimental data is comparable to the correlation of the FE reference solution.\\
The failure loads predicted for the humeri FFH5L and FFH6R were reproduced with very good agreement, and are of similar accuracy as the failure load of the calibrated bone FFH5R with relative error below $3.8\:\%$ for all bones. In \cite{dahan2020neck}, a classic maximum principal compression strain criteria did not reproduce well the experimental yield loads. The phase-field FCM model presents itself as a promising alternative. The comparison of load-displacement curves shows, that the proposed model is able to capture the non-brittle failure behavior observed experimentally, even though it has originally been designed for brittle fracture. For all three humeri, there is a visible "plastic" phase between the onset of damage and the failure point. The proposed PFM is able to capture this behavior due to the choice of the degradation function and the length-scale parameter. A quadratic degradation function allows for an accumulation of damage before the onset of fracture \cite{sargado2018high}, which is undesired when reproducing brittle fracture but allows to mimic a "plastic" damage phase representing a more ductile behavior. Choosing a larger $l_0$ leads to an increased region of damage around the crack which in turn reinforces the premature decrease in material stiffness. Shen et al.~\cite{shen2019novel} suggested to use a generic degradation function, which interpolates between a quadratic and a cubic degradation function. Based on an parameter $\bar{s}$ the quadratic degradation function with a more ductile behavior is recovered for $\bar{s}\,=\,2$, while a brittle behavior can be obtained for $\bar{s} \rightarrow 0$. They suggest a choice of $\bar{s}\,=\,0.8$ which represents a quasi-brittle response. Moreover, they obtain a length-scale parameter of $l_{0_\textrm{Shen}}=1.5$ when calibrating their model which is similar to $l_0=2.25$ found in our work. Interestingly, the experiment used for calibration in Shen et al.~\cite{shen2019novel} involves a fracture at the anatomical neck of a humerus. As stated in Dahan \cite{dahan2020neck}, anatomical neck fractures show a more brittle behavior than the surgical neck fractures considered in this contribution. This fact can very well explain the different choices in the degradation function and the length-scale parameter.\\
A qualitative comparison of fracture paths confirms, that the PFM can reproduce the surgical neck fractures as observed experimentally. In agreement with the load-displacement curves, damage accumulates in the trabecular regions inside the humeral head. However, a clear initiation of the fracture inside the humeral head and not on the outer cortex as suspected by Dahan~\cite{dahan2020neck} can only be observed numerically for the humeri FFH6R. A larger number of humeri would be necessary to make justified statements about the initiation process, including advanced experimental techniques to analyze the point of initiation not only numerically but experimentally. Including more experiments could decisively strengthen the validity of the study. Moreover, there are some uncertainties in the material parameters used. Even though the dependency of the elastic strain energy density $G_c$ on the Young's modulus has shown to be a valid choice in this setting, its distribution in the bone and its exact influence on fracture patterns and failure loads is still an open question. \textcolor{changes}{In addition, the choice of the volumetric-deviatoric split by Amor~et~al.~\cite{amor2009regularized} for the PFM under compressive load may be questioned. Following Eq.~\ref{eq:governeq}a, only the positive part of the elastic strain energy is degraded and contributes to crack propagation. However, an analysis of principle stresses along the crack surface (Appendix~\ref{appendix:a}) confirms that tensile stresses are present that drive the fracture, and the results in Section \ref{sec:results} demonstrate that the fracture behavior can be captured. A detailed investigation and comparison of the different tension-compression splits is beyond the scope of the present contribution and subject of future research.}

\section{Conclusion} \label{sec:conclusion}
 \textcolor{changes}{In this contribution, the prediction of proximal humerus fractures with a PFM is investigated. The proposed numerical framework combines a PFM with a spatially varying critical energy release rate~\cite{shen2019novel} and an embedded domain approach for an efficient handling of the complex bone geometry. A validation on three in vitro experiments inducing surgical neck fractures in human humeri~\cite{dahan2020neck} is presented. To the authors' knowledge, this is the first validation for neck fractures in human humeri including a comparison of DIC strains on the bones' surface, experimental failure loads, and a quantitative comparison of the fracture paths.} \\
After calibration of the phase-field length parameter based on one humeri, strains in the linear elastic regime, failure loads and fracture paths are compared for the three bone specimen. Strains on the bones' surface agree moderately well with experimental data and are of similar accuracy compared to a FE reference solution ~\cite{dahan2020neck}. The proposed model is capable of reproducing the non-brittle behavior in the experimental load-displacement curves which materializes in a plastic phase after the onset of damage. A qualitative comparison of fracture patterns shows that for all three humeri failures occur at the surgical neck of the humeral head. Moreover, the fracture paths agree well with experiments. Additionally, the computed failure loads show an excellent agreement with experimental data with \textcolor{changes}{a relative error below $3.8\,\%$} for all three humeri. The results demonstrate that phase-field models present a promising tool in the patient-specific prediction of surgical neck fractures in human humeri.

%% file: utils/acknowledgements.tex
LH, SK and ER gratefully acknowledge the funding through Deutsche Forschungsgemeinschaft (DFG) for its financial support through the TUM International Graduate School of Science and Engineering (IGSSE), GSC 81.

%% file: utils/appendix.tex
\appendix
\cleardoublepage\phantomsection\addcontentsline{toc}{section}{Appendix}
\section*{Appendix}
\section{Principle Stress Analysis}\label{appendix:a}
\begin{figure}[h!]
	\centering
	\includegraphics[width=0.95\textwidth]{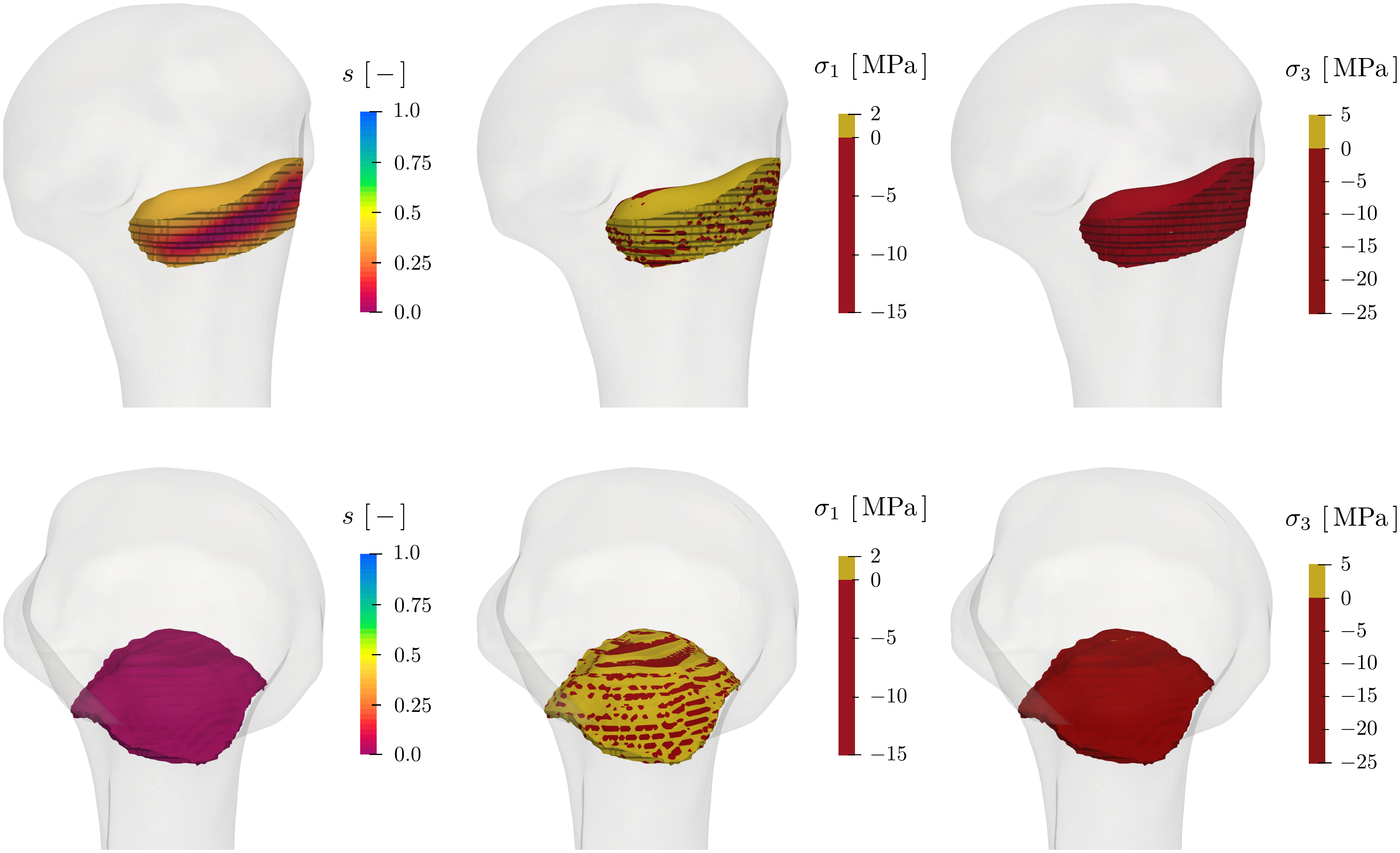}
	\caption{The phase-field $s$, the maximum and the minimum principal stresses $\sigma_1$ and $\sigma_3$ are visualized for the humeri FFH5R. In the top row, the point of crack initiation at a displacement of $u_z=-0.4675$ mm is depicted on an iso-volume of the phase-field with $0\leq s \leq 0.4$. In the bottom row, the fully developed crack surface is shown using an iso-volume of the phase-field with $0\leq s \leq 0.01$.}
	\label{fig:stressanalysis}
\end{figure}

To justify the choice of a phase-field formulation based on a full tension-compression split with a volumetric-deviatoric decomposition~\cite{amor2009regularized} an analysis of the principal stresses is presented in Figure~\ref{fig:stressanalysis}. The phase-field $s$ along with the maximum and minimum principal stresses $\sigma_1$ and $\sigma_3$ are plotted for the crack initiation (top) and on the full crack surface (bottom) for the humeri FFH5R. 
